\documentclass[preprint,3p,12pt,pdf]{elsarticle}

\usepackage{lineno,hyperref}
\modulolinenumbers[5]

\usepackage{hyperref}

\usepackage{epstopdf}

\usepackage{tikz}
\usetikzlibrary{arrows}

\usepackage{pst-node}
\usepackage{tikz-cd}

\usepackage[shellescape]{gmp}

\usepackage{mathrsfs}
\usepackage{amssymb}
\usepackage{amsfonts}
\usepackage{latexsym}
\usepackage{mathtools}
\usepackage{xcolor}
\usepackage{wrapfig}
\usepackage{floatflt}

\usepackage{mathtools}
\usepackage{extarrows}

\usepackage{graphicx}

\usepackage{subcaption}

\def\lb{\label}

\newcommand{\er}[1]{\textrm{(\ref{#1})}}

\def\a{\alpha}         
\def\b{\beta}          
         
\def\G{\Gamma}

\def\l{\lambda}        
\def\L{\Lambda}

\def\r{\rho}

\def\o{\omega}

       \def\vp{\varphi}    

\def\Z{{\mathbb Z}}    \def\R{{\mathbb R}}   \def\C{{\mathbb C}}    
    \def\N{{\mathbb N}}   

\def\lt{\biggl}                  \def\rt{\biggr}
\def\ol{\overline}               \def\wt{\widetilde}

\let\ge\geqslant                 \let\le\leqslant

\def\iy{\infty}

\def\pa{\partial}

\def\el2{\ell^{\,2}}             \def\1{1\!\!1}

\def\Re{\mathop{\mathrm{Re}}\nolimits}

\let\ge\geqslant
\let\le\leqslant

\newcommand{\ca}{\begin{cases}}
	\newcommand{\ac}{\end{cases}}
\newcommand{\ma}{\begin{pmatrix}}
	\newcommand{\am}{\end{pmatrix}}
\def\eq{\begin{equation}}
	\def\qe{\end{equation}}
\def\[{\begin{equation}}
	\def\]{\end{equation}}

\begin{document}
	
	\begin{frontmatter}

		\title{Periodic oscillations of coefficients of
			power series that satisfy functional equations, a practical revision}

		\date{\today}

		\author
		{Anton A. Kutsenko}
	
		\address{Mathematical Institute for Machine Learning and Data
Science, KU Eichst\"att--Ingolstadt, Germany; email: akucenko@gmail.com}
	
	\begin{abstract}
		For the solutions $\Phi(z)$ of functional equations $\Phi(z)=P(z)+\Phi(Q(z))$, we derive a complete asymptotic of power series coefficients. As an application, we improve significantly an asymptotic of the number of $2,3$-trees with $n$ leaves given in Adv. Math. 44:180–205, 1982 by Andrew M. Odlyzko. The methods we consider can be applied to more general functional equations too.
	\end{abstract}


	
\end{frontmatter}


{\section{Introduction}\lb{sec0}}

In the good classic work \cite{O} of 1982, the author considers the functional equation
\[\lb{int0}
 f(z)=z+f(z^2+z^3)
\]
in the context of finding asymptotics of the number of $2,3$-trees with $n$ leaves. The author proves that the Taylor coefficients $a_n$ of $f(z)$, which represent these numbers of trees, satisfy the asymptotic
\[\lb{int1}
 a_n\sim \frac{\vp^n}{n}u(\ln n),
\]
where $\vp$ is the golden ratio, and $u$ is some analytic $\ln(4-\vp)$-periodic function with the mean $(\vp\ln(4-\vp))^{-1}$. Then, the author of \cite{O} extends such type of results to more general functional equations \er{001}. However, the author of \cite{O} note `` unfortunately we do not obtain any good expansions for u(x)". This phrase and the next one:

``
The proof that is presented here yields a result somewhat stronger than the assertion of the theorem, namely, that 
$$
 a_n=\frac{\vp^n}{n}u(\ln n)+O\lt(\frac{\vp^n}{n^2}\rt),\ \ {\rm as}\ \ n\to\iy.
$$ 
With additional work one can obtain an even more complete asymptotic 
development of the $a_n$. A problem which is left open by the present proof of Theorem 2 is that of obtaining an explicit representation of the periodic function $u(x)$. 
"

is the main motivation for the current research. In this current work, we provide a complete asymptotic of Taylor coefficients of the solution of the general functional equation \er{001}, see the results in \er{029}-\er{033}. Perhaps, these results admit further simplification, but, it is enough for our current needs. All the functions involved in the results admit straightforward numerical implementation as demonstrated in Section \ref{sec2}. Some of the methods presented in this work have some similarities with those proposed in the recently published paper \cite{Ku} devoted to the Schr\"oder-type functional equations. In fact, the methods can be applied to slightly more general functional
equations, see {\bf Remark} at the end of Section \ref{sec1}. Finally note that there are many very nice works devoted to the topic, some of them are mentioned in \cite{Ku} and \cite{T}, see also papers which cite \cite{T}, and we do not refer to them again in this version of the paper. Most of these works focus on the first or few first terms in asymptotics of power series coefficients considering sometimes more general functional equations than \er{001}. We focus mostly on \er{001}, but our goal is the complete asymptotic series.

{\section{Main results}\lb{sec1}}

 Consder the functional equation
\[\lb{001}
 \Phi(z)=P(z)+\Phi(Q(z)),
\]
where $P$ and $Q$ are some non-zero analytic functions, defined on a sufficiently large neighborhood of $0$. The functions $P$ and $Q$ are assumed to be real-valued for real arguments. The most practical case is when $P$ and $Q$ are entire functions, e.g. polynomials. For the assumptions regarding $P$ and $Q$, we follow that discussed in \cite{O}. Generally speaking, one of the key assumption is that the filled Julia set for $Q(z)$ contains the ball $\{|z|\le q\}$, and the $Q$-iterations of all the points in this ball tend to $0$ except for this unique point $q>0$. In fact, the last condition can be weakened but we consider the simplest case. There are only two fixed points $Q(0)=0$ and $Q(q)=q$ in the ball - the first one is an attracting point with $0\le Q'(0)<1$ and the second one is a repelling point with $Q'(q)>1$. The assumptions on $P$ are less strict, some of them are $P(0)=0$ and $P(q)\ne0$. Other conditions will appear as needed. We have 
$$
 Q_N(z):=\underbrace{Q\circ...\circ Q}_{N}(z)\to0
$$ 
for any $z$ such that $|z|\le|q|$ and $z\ne q$. Then, the solution of \er{001} can be written as
\[\lb{002}
 \Phi(z)=P(z)+P(Q(z))+P(Q_2(z))+P(Q_3(z))+....
\]
For $|z|<q$, this function can be expanded into the Taylor series
\[\lb{003}
 \Phi(z)=\vp_1z+\vp_2z^2+\vp_3z^3+....
\]
The goal is to determine an explicit asymptotic of $\vp_n$ for $n\to\iy$. The first step is to find $T(z)$ satisfying the functional equation and the initial conditions
\[\lb{004}
 T(z)=e^{\a P(z)}T(Q(z)),\ \ \ T(q)=0,\ \ \ T'(q)=-1.
\]
Differentiating $T(z)$ in \er{004}, we obtain
\[\lb{005}
 \a=-\frac{\ln Q'(q)}{P(q)}.
\]
There exists a solution of \er{004} analytic in some neighborhood of $q$, logarithm of which can be expanded into the series
\[\lb{006}
 \ln T(z)=\ln(q-z)+\sum_{m=1}^{+\iy}\lt(\ln\frac{q-Q_{-m}(z)}{q-Q_{1-m}(z)}-\a P(Q_{-m}(z))\rt),
\] 
where
$$
 Q_{-N}(z):=\underbrace{Q^{-1}\circ...\circ Q^{-1}}_{N}(z)\to q
$$ 
for $z$ in some neighborhood of $q$, since $Q'(q)>1$, and, hence $q$ is an attracting point for $Q^{-1}$. Using, e.g., the Newton method, the inverse function $Q^{-1}$ can be computed numerically. We assume also that this function exists. The series
\er{006} can be rewritten in the form
\[\lb{007}
 \ln T(z)=\ln(q-z)+\sum_{m=1}^{+\iy}\ln\frac{Q'(q)^{\frac{P(Q_{-m}(z))}{P(q)}}}{R(Q_{-m}(z))},
\]
where the analytic function $R$ is defined by
\[\lb{008}
 R(z):=\frac{Q(z)-q}{z-q},\ \ \ R(q):=Q'(q).
\]
Series \er{006} and \er{008} converge exponentially fast, and their numerical implementation is straightforward. The function
\[\lb{009}
 \wt\L(z):=\Phi(z)-\frac1{\a}\ln T(z)
\]
satisfies the functional equation
\[\lb{010}
 \wt \L(Q(z))=\wt\L(z),
\]
see \er{001} and \er{004}. Let us define the analytic function
\[\lb{011}
 \Pi(z):=\lim_{N\to+\iy}Q_N\lt(q-\frac{z}{Q'(q)^N}\rt),
\]
satisfying the Poincar\'e-type functional equation and initial conditions
\[\lb{012}
 \Pi(Q'(q)z)=Q(\Pi(z)),\ \ \ \Pi(0)=q,\ \ \ \Pi'(0)=-1.
\]
If $Q$ is entire then $\Pi$ is entire as well. For the numerical computation of $\Pi$ one may use the recurrence scheme based on the identity
$$
 Q_N\lt(q-\frac{1}{Q'(q)^N}\cdot z\rt)=Q_{N-1}\lt(q-\frac{1}{Q'(q)^{N-1}}\cdot\frac{zR(q-\frac{z}{Q'(q)^N})}{Q'(q)}\rt),
$$
see \er{008}. The corresponding fast convergent recursion can be programmed relatively easily. Combining \er{010} and \er{012}, we obtain the $1$-periodic function
\[\lb{013}
 \L(z):=\wt\L(\Pi(Q'(q)^z)),
\]
which can be expanded into the Fourier series
\[\lb{014}
 \L(z)=\sum_{m=-\iy}^{+\iy}\l_me^{2\pi\mathbf{i}mz},
\]
where $\mathbf{i}=\sqrt{-1}$ and $\l_{-m}=\ol{\l_m}$. Without providing the details, we note that $\L(z)$ is analytic in some strip symmetric about the real axis. The width of the strip depends on the geometric properties of the Julia set related to the function $Q(z)$. The most significant influence on the width of the strip is the geometric structure of the Julia set in the vicinity of $q$.  Using \er{009} along with\er{013}, and \er{014} we deduce that
\[\lb{015}
 \Phi(z)=\L(\frac{\ln\Psi(z)}{\ln Q'(q)})+\frac1{\a}\ln T(z)=\frac1{\a}\ln T(z)+\sum_{m=-\iy}^{+\iy}\l_m\Psi(z)^{\frac{2\pi\mathbf{i}m}{\ln Q'(q)}},
\]
where the analytic function
\[\lb{016}
 \Psi(z):=\Pi^{-1}(z)=\lim_{N\to+\iy}Q'(q)^N(q-Q_{-N}(z))
\]
satisfies the Schr\"oder-type functional equation
\[\lb{017}
 \Psi(Q(z))=Q'(q)\Psi(z),\ \ \ \Psi(q)=0,\ \ \ \Psi'(q)=-1.
\]
The Taylor series of this function is
\[\lb{018}
 \Psi(z)=(q-z)+\frac{\Psi''(q)}{2!}(q-z)^2-\frac{\Psi'''(q)}{3!}(q-z)^3+...,
\]
where the derivatives can be found by the differentiation of \er{017}, namely $\Psi''(q)=\frac{-Q''(q)}{Q'(q)-Q'(q)^2}$ and, generally, the Fa\`a di Bruno's formula applied to \er{017} gives the recurrence identity
\[\lb{019}
 \Psi^{(m)}(q)=\frac1{Q'(q)-Q'(q)^m}\sum_{k=1}^{m-1}\Psi^{(k)}(q)B_{m,k}(Q'(q),Q''(q),...,Q^{(m-k+1)}(q))
\]
with $B_{m,k}$ denoting the Bell polynomials
\[\lb{Bell1}
 B_{m,k}(x_1,x_2,...,x_{m-k+1})=\sum\frac{m!}{j_1!j_2!...j_{m-k+1}!}\lt(\frac{x_1}{1!}\rt)^{j_1}\lt(\frac{x_2}{2!}\rt)^{j_2}...\lt(\frac{x_{m-k+1}}{(m-k+1)!}\rt)^{j_{m-k+1}},
\]
where the sum for Bell polynomials is taken over all sequences $j_r$ of non-negative integers such that
\[\lb{Bell2}
 j_1+j_2+...+j_{m-k+1}=k,\ \ \ j_1+2j_2+3j_3+...+(m-k+1)j_{m-k+1}=m.
\]
There are some examples of Bell polynomials
\begin{multline}\lb{Bell3}
 B_{m,1}(x_1,...,x_m)=x_m,\ \ \ B_{m,2}(x_1,...,x_{m-1})=\frac12\sum_{j=1}^{m-1}\binom{m}{j}x_jx_{m-j},\\ B_{m,m-1}(x_1,x_2)=\binom{m}{2}x_1^{m-2}x_2,\ \ \ B_{m,m}(x_1)=x_1,
\end{multline}
where $\binom{a}{b}$ denotes the standard binomial coefficients. Applying the Fa\`a di Bruno's formula to \er{018}, we obtain also
\[\lb{020}
 \Psi(z)^r=(q-z)^r+\frac{\psi_1(r)}{1!}(q-z)^{r+1}+\frac{\psi_2(r)}{2!}(q-z)^{r+2}+...,
\]
where $\psi_1(r)=\frac{-rQ''(q)}{2(Q'(q)-Q'(q)^2)}$ and, generally, the polynomials $\psi_m(r)$ are given by
\[\lb{021}
 \psi_m(r)=\sum_{k=1}^{m}\binom{r}{k}B_{m,k}\lt(-\frac{\Psi''(q)}{2},-\frac{\Psi'''(q)}{3},...,-\frac{\Psi^{(m-k+2)}(q)}{m-k+2}\rt),
\]
with $\binom{r}{k}$ denoting the generalized binomial coefficients 
\[\lb{022}
 \binom{r}{k}=\frac{r(r-1)...(r-k+1)}{k!},\ \ \ r\in\C.
\]
Substituting \er{020} into \er{015} and denoting $\psi_0(r)=1$ and $\b=\ln Q'(q)$, we obtain
\begin{multline}\lb{023}
 \Phi(z)=\frac1{\a}\ln T(z)+\sum_{m=-\iy}^{+\iy}\l_m\sum_{j=0}^{+\iy}\frac{\psi_j(\frac{2\pi\mathbf{i}m}{\b})}{j!}(q-z)^{j+\frac{2\pi\mathbf{i}m}{\b}}=\\
 =\frac1{\a}\ln T(z)+\sum_{m=-\iy}^{+\iy}\l_m\sum_{j=0}^{+\iy}\frac{\psi_j(\frac{2\pi\mathbf{i}m}{\b})}{j!}q^{j+\frac{2\pi\mathbf{i}m}{\b}}\sum_{n=0}^{+\iy}(-q)^{-n}\binom{j+\frac{2\pi\mathbf{i}m}{\b}}{n}z^n.
\end{multline}
Assuming that the Fourier coefficients $\l_m$ tend to zero sufficiently fast, and taking into account the fact that $q$ is a unique singularity of $\Phi$ in the ball $|z|\le q$, we can equate the coefficients with the same $z^n$ in \er{023} to obtain the asymptotic of $\vp_n$. The $\sum$ in \er{007} can be dropped out, since it is analytic in the neighborhood of $q$ and the Taylor coefficients of this $\sum$ have exponential growth (or attenuation) less than $q^n$. Thus, we have
\[\lb{024}
 \vp_n\simeq\frac{-1}{\a nq^n}+\sum_{m=-\iy}^{+\iy}\l_m\sum_{j=0}^{+\iy}\frac{\psi_j(\frac{2\pi\mathbf{i}m}{\b})}{j!}q^{j+\frac{2\pi\mathbf{i}m}{\b}}(-q)^{-n}\binom{j+\frac{2\pi\mathbf{i}m}{\b}}{n}.
\]
For the binomial coefficients, there exists an asymptotic formula
\[\lb{025}
 (-1)^n\binom{r}{n}\simeq\frac{1}{\Gamma(-r)n^{r+1}}\lt(\frac{S_{0}(r)}{n^0}+\frac{S_{2}(r)}{n^1}+\frac{S_{4}(r)}{n^2}+\frac{S_{6}(r)}{n^3}...\rt),
\]
where
\[\lb{026}
 S_0(r)=1,\ \ \ S_2(r)=\frac{r(r+1)}2,\ \ \ S_4(r)=\frac{r(r+1)(r+2)(3r+1)}{24},
\]
and, generally, the polynomials $S_{2k}(r)$ of degree $2k$ can be defined from the identity $\binom{r}{n+1}=\frac{r-n}{n+1}\binom{r}{n}$, which gives
\begin{multline}\lb{027}
 (k-1)S_{2(k-1)}(r)=S_{2(k-2)}(r)\lt(\binom{-(k-2)-r-1}{2}-(r+1)\rt)+S_{2(k-3)}(r)\cdot\\
 \cdot\lt(\binom{-(k-3)-r-1}{3}+(r+1)\rt)+...+S_{0}(r)\lt(\binom{-r-1}{k}+(-1)^{k-1}(r+1)\rt).
\end{multline}
Substituting \er{025} into \er{024} and using standard properties of $\G$-function, we obtain
\begin{multline}\lb{028}
 q^n\vp_n\simeq\frac{-1}{\a n}+\sum_{m=-\iy}^{+\iy}\l_m\sum_{j=0}^{+\iy}\frac{\psi_j(\frac{2\pi\mathbf{i}m}{\b})}{j!}\frac{q^{j+\frac{2\pi\mathbf{i}m}{\b}}}{\Gamma(-j-\frac{2\pi\mathbf{i}m}{\b})n^{1+j+\frac{2\pi\mathbf{i}m}{\b}}}\sum_{k=0}^{+\iy}\frac{S_{2k}(j+\frac{2\pi\mathbf{i}m}{\b})}{n^k}=\\
 =\frac{-1}{\a n}+\sum_{m=-\iy}^{+\iy}\frac{\l_m}{\Gamma(-\frac{2\pi\mathbf{i}m}{\b})}\sum_{j=0}^{+\iy}\frac{\psi_j(\frac{2\pi\mathbf{i}m}{\b})\binom{-1-\frac{2\pi\mathbf{i}m}{\b}}{j}q^{j+\frac{2\pi\mathbf{i}m}{\b}}}{n^{1+j+\frac{2\pi\mathbf{i}m}{\b}}}\sum_{k=0}^{+\iy}\frac{S_{2k}(j+\frac{2\pi\mathbf{i}m}{\b})}{n^k}=\\
 =\frac{-1}{\a n}+\sum_{r=0}^{+\iy}\frac1{n^{1+r}}\sum_{m=-\iy}^{+\iy}\frac{\l_me^{2\pi\mathbf{i}m\frac{\ln q-\ln n}{\b}}}{\Gamma(-\frac{2\pi\mathbf{i}m}{\b})}
 \sum_{\substack{j,k\ge0\\ {j+k=r}}}q^{j}\psi_j(\frac{2\pi\mathbf{i}m}{\b})\binom{\frac{-\b-2\pi\mathbf{i}m}{\b}}{j}S_{2k}(\frac{j\b+2\pi\mathbf{i}m}{\b}).
\end{multline}
Asymptotic expansion \er{028} is the main result. Using \er{005} and some other formulas presented above, it is convenient to rewrite \er{028} in the form
\[\lb{029}
 \vp_n\simeq\frac{q^{-n}}{n}K_1(\frac{\ln q-\ln n}{\b})+\frac{q^{-n}}{n^2}K_2(\frac{\ln q-\ln n}{\b})+\frac{q^{-n}}{n^3}K_3(\frac{\ln q-\ln n}{\b})+...,
\]
where $1$-periodic real functions $K_r(x)$, $x\in\R$ are given by
\[\lb{030}
 K_1(x)=-\frac{1}{\a}+2\Re\sum_{m=1}^{+\iy}\frac{\l_me^{2\pi\mathbf{i}mx}}{\Gamma(-\frac{2\pi\mathbf{i}m}{\b})},
\]
\[\lb{031}
 K_2(x)=2\Re\sum_{m=1}^{+\iy}\frac{\l_me^{2\pi\mathbf{i}mx}}{\Gamma(-\frac{2\pi\mathbf{i}m}{\b})}\cdot\frac{2\pi\mathbf{i}m(2\pi\mathbf{i}m+\b)(qQ''(q)+Q'(q)-Q'(q)^2)}{2\b^2(Q'(q)-Q'(q)^2)},
\]
and, generally, 
\[\lb{032}
 K_r(x)=2\Re\sum_{m=1}^{+\iy}\frac{\l_me^{2\pi\mathbf{i}mx}}{\Gamma(-\frac{2\pi\mathbf{i}m}{\b})}A_r(2\pi\mathbf{i}m),\ \ \ r\ge2,
\]
where the polynomials $A_r(z)$ are given by
\[\lb{033}
 A_r(z)=\sum_{\substack{j,k\ge0\\ {j+k=r}}}q^{j}\psi_j(\frac{z}{\b})\binom{\frac{-\b-z}{\b}}{j}S_{2k}(\frac{j\b+z}{\b}).
\]
Noting that the derivative $\pa/\pa x$ leads to the multiplication by $2\pi\mathbf{i}m$ of $m$-th term in the Fourier series, one may express $K_r(x)$ through the $A_r(\pa/\pa x)K_1(x)$. 

Everything is ready for the numerical implementation. However, it should be noted that the Fourier coefficients $\l_m$ decay exponentially fast, and $\Gamma(-\frac{2\pi\mathbf{i}m}{\b})$ also decay exponentially fast. Thus, we have the ratio of two small quantities in \er{030}-\er{032}. To improve the computation of this ratio, one may compute Fourier coefficients $\hat \l_m=e^{2\pi m y}\l_m$ of $\L(x-\mathbf{i}y)$, where $y>0$ and $m\ge0$. A good strategy is to find the maximal $y$ for which the computations still give proper results, without extremely large, discontinuous, or NaN values. The parameter $y$ should be greater than $\frac{\pi}{2\b}$, since
$$
 |\Gamma(\mathbf{i}x)|^2=\frac{\pi}{x\sinh\pi x},\ \ \ x\in\R.
$$
The ratio of two small quantities can be rewritten as
\[\lb{034}
 \frac{\l_m}{\Gamma(-\frac{2\pi\mathbf{i}m}{\b})}=\exp(\ln\hat\l_m+(\ln\Gamma(-\frac{2\pi\mathbf{i}m}{\b})-2\pi my)),
\]
which gives very accurate results as tested in numerical examples. Not also that the computation of Fourier coefficients itself can be done by applying FFT to the array of $\L(x-\mathbf{i}y)$ with $x=-n,-n+1/N,...,-n+(N-1)/N$ with some $n\in\N$ and large $N\in\N$. Moreover, all the procedures described above admit a straightforward vectorization to compute the array $\L(x-\mathbf{i}y)$ quickly.

{\bf Remark.} In fact, there are many papers devoted to the first or first few asymptotic terms of power series coefficients of solutions of various functional equations, see, e.g. the reference list in \cite{T}, and the papers which cite \cite{T}. The key point of our current research is to obtain a complete asymptotic series, see \er{029}. This is a reason why we focus on equations of the certain type \er{001}. However, the methods we use are applicable to, e.g., a slightly more general functional
equation
$$
 f(z)=m(z)f(\vp(z))+a(z),
$$
because the corresponding power series coefficients can be expressed as
$$
 f_n=\sum_{k\in\Z}\frac{\hat\o(k)\r^{2k\pi\mathbf{i}\kappa+\a}}{(\eta\rho)^n}\binom{2k\pi\mathbf{i}\kappa+\a}{n}
$$
with some $\hat\o$, $\r$, $\kappa$, $\a$, and $\eta$, see Proof of Theorem 19 in \cite{T}. Now, using asymptotic series \er{025} for the binomial coefficients one can arrive to the asymptotic series of $f_n$ similar to \er{029}. This is the most important point. Obtaining explicit formulas for $\hat\o$, $\r$, $\kappa$, $\a$, and $\eta$ of the same form as in \er{019}, \er{021} and etc. also seems to be quite realizable.

{\section{Example}\lb{sec2}}

We apply the results \er{029}-\er{033} to the equation \er{int0}. In this case
\[\lb{100}
 Q(z)=z^2+z^3,\ \ \ P(z)=z,\ \ \ q=\vp^{-1}=\frac{\sqrt{5}-1}2,\ \ \ \b=\ln(4-\vp),\ \ \ \a=-\vp\ln(4-\vp).
\]
The normalized Fourier coefficients $\hat \l_m=e^{4\pi m}\l_m$, see \er{014} and the discussion at the end of Section \ref{sec1}, for $m=1,...,10$ are
$$
 -0.10417+0.0052295\mathbf{i},\ \    
 0.10883+0.04913\mathbf{i},\ \ 
 -0.0027473+0.02632\mathbf{i},\ \ 
  0.011381+0.0076878\mathbf{i}, 
 $$
$$ 
-0.0010885+0.0032545\mathbf{i},\ \
 0.00099529+0.0001076\mathbf{i},\ \
 -0.0013305-0.00023601\mathbf{i}, 
$$
$$  
-0.00059214+0.00054537\mathbf{i},\ \
 0.00007277+0.00032196\mathbf{i},\ \
 0.000088894-0.000088991\mathbf{i}.
$$
The corresponding ratios with $\G$-function computed by \er{034} with $y=-2$ are respectively
$$
   -0.033869+0.0013274\mathbf{i},\ \ 
   0.0047334-0.015924\mathbf{i},\ \ 
   -0.00061251+0.0012199\mathbf{i}, 
$$
$$
    0.00017226-0.00017793\mathbf{i},\ \
    0.000017638+0.000011296i\mathbf{i},\ \
   0.0000019278+0.00000062387i\mathbf{i},
$$
$$   
    8.9\cdot10^{-7}-1.7\cdot10^{-8}\mathbf{i},\ \
   -1.6\cdot10^{-7}+5.1\cdot10^{-8}\mathbf{i},\ \
   2.2\cdot10^{-8}+4.4\cdot10^{-9}\mathbf{i},\ \
   -2.5\cdot10^{-9}-1.1\cdot10^{-9}\mathbf{i}.
$$
They decay relatively fast. We use these values for the computation of $1$-periodic functions $K_1(x)$ and $K_2(x)$, see Fig. \ref{fig1}.

\begin{figure}[h]
	\center{\includegraphics[width=0.9\linewidth]{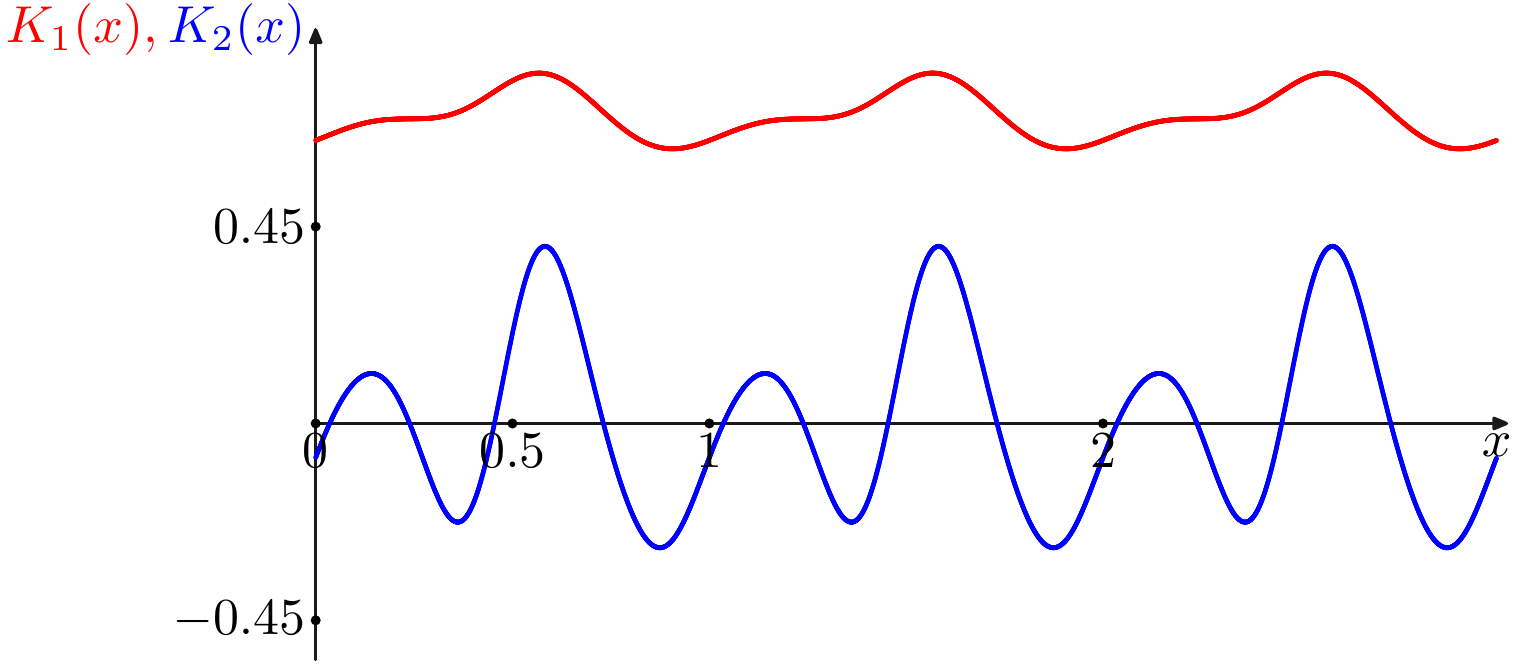}}
	\caption{Periodic functions $K_1(x)$ and $K_2(x)$, see \er{030} and \er{031}, computed for the case \er{100}.}\lb{fig1}
\end{figure}

The exact values $\vp_n$ or their re-scaled analogs $n\vp^{-n}\vp_n$ can be computed by using the recurrence identity
$$
 \vp_n=\sum_{2k+3m=n}\binom{k+m}{k}\vp_{k+m},
$$ 
which follows directly from \er{int0}. The first $10000$ values $\vp_n$ are compared with asymptotic terms in Fig. \ref{fig2}. Asymptotic terms fit exact values very well.  

Modules of highly optimized code with active use of parallel programming for computing the functions described in Sections \ref{sec1} and \ref{sec2} are available on a special request to AK.

\begin{figure}[h]
	\center{\includegraphics[width=0.9\linewidth]{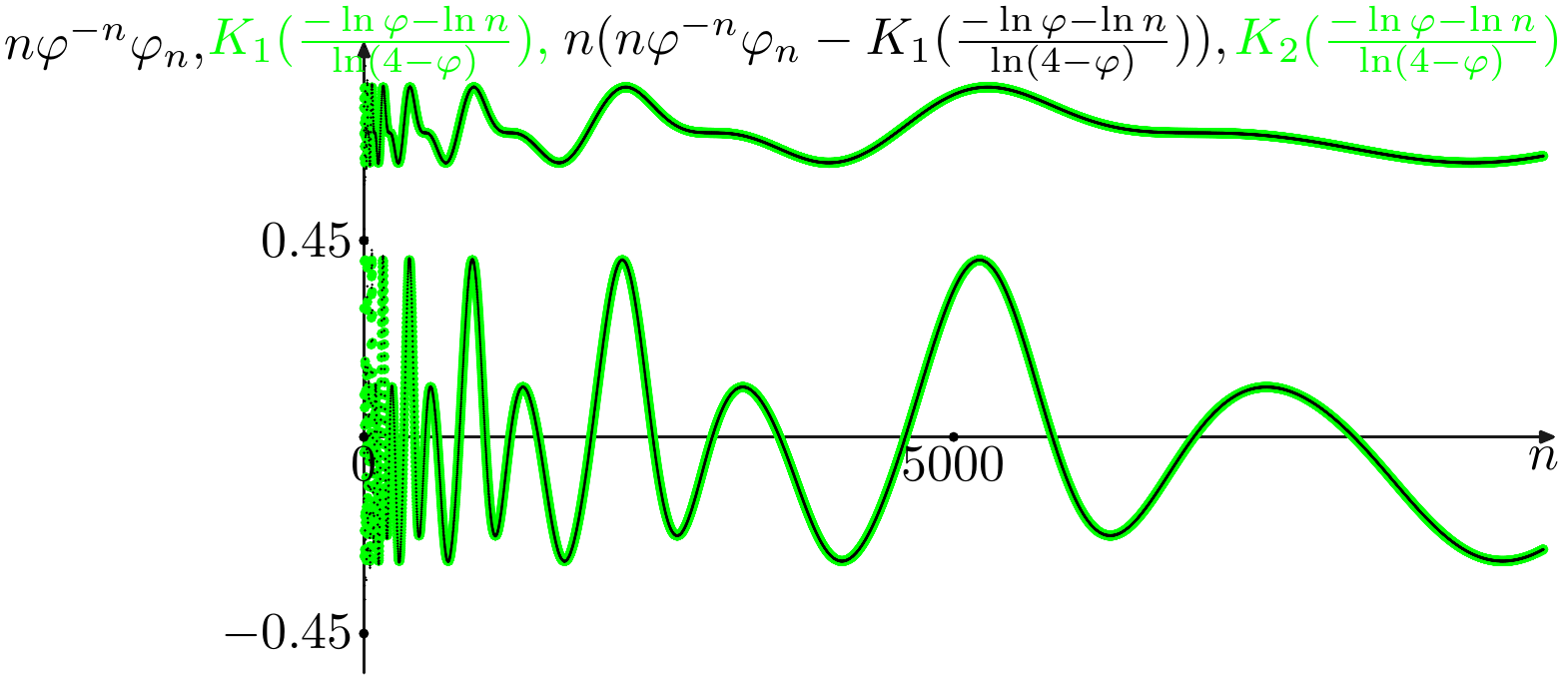}}
	\caption{Exact normalized values $\vp_n$, and two first asymptotic terms in \er{029}.
	}\lb{fig2}
\end{figure}

\section*{Acknowledgements} 
This paper is a contribution to the project M3 of the Collaborative Research Centre TRR 181 "Energy Transfer in Atmosphere and Ocean" funded by the Deutsche Forschungsgemeinschaft (DFG, German Research Foundation) - Projektnummer 274762653. 

\end{document}